\documentclass[12pt]{amsart}
\usepackage{enumerate,amstext,amsfonts,amssymb,amscd,amsbsy,amsmath,verbatim}
\usepackage{ifthen,multicol,extarrows,fullpage,anyfontsize,mathtools,multirow}
\usepackage[lite,alphabetic]{amsrefs} 
\usepackage{lscape}
\usepackage{enumitem}
\usepackage{color,tikz}
\usepackage{amsthm}
\usepackage{latexsym}
\usepackage[all]{xy}
\usepackage[bookmarks, colorlinks=true, linkcolor=blue, citecolor=blue, urlcolor=blue]{hyperref}
\usepackage{tikz-cd}

\usepackage{eucal}

\numberwithin{equation}{section}

\newtheorem{claim*}{Claim}

\theoremstyle{remark}

\newcommand{\Tor}{\operatorname{Tor}}

\newcommand{\Sym}{\operatorname{Sym}} 
\newcommand{\GL}{\mathbf{GL}}

\newcommand{\defi}[1]{\textsf{#1}} 

\newcommand{\PP}{\mathbb P}

\newcommand{\ZZ}{\mathbb Z}

\newcommand{\CC}{\mathbb C}
\newcommand{\ba}{\mathbf a}

\newcommand{\bS}{\mathbf{S}}

\newcommand{\cO}{\mathcal O}

\usepackage{rotating}

\title{The Schur-Veronese package in Macaulay2}
\author{Juliette Bruce}
\address{Department of Mathematics, University of California, Berkeley, Berkeley, CA}
\email{\href{mailto:juliette.bruce@berkeley.edu}{juliette.bruce@berkeley.edu}}
\urladdr{\url{https://juliettebruce.github.io}}

\author{Daniel Erman}
\address{Department of Mathematics, University of Wisconsin, Madison, WI}
\email{\href{mailto:derman@math.wisc.edu}{derman@math.wisc.edu}}
\urladdr{\url{http://math.wisc.edu/~derman/}}

\author{Steve Goldstein}
\address{Botany Department and Department of Biostatistics and Medical Informatics, University of Wisconsin, Madison, WI}
\email{\href{mailto:sgoldstein@wisc.edu}{sgoldstein@wisc.edu}}

\author{Jay Yang}
\address{School of Mathematics, University of Minnesota, Minneapolis, MN}
\email{\href{mailto:jkyang@umn.edu}{jkyang@umn.edu}}
\urladdr{\url{http://www-users.math.umn.edu/~jkyang/}}
\date{\today}


\begin{document}
\thanks{
JB received support from the NSF GRFP under grant DGE-1256259, NSF grant DMS-1502553, NSF MSPRF DMS-2002239, and from the Graduate School and the Office of the Vice Chancellor for Research and Graduate Education at the University of Wisconsin-Madison with funding from the Wisconsin Alumni Research Foundation. JB is grateful for the support of the Mathematical Sciences Research Institute in Berkeley, California, where she was in residence for the Fall 2020 semester. DE received support from NSF grant DMS-1601619. JY received support from NSF grant DMS-1502553.
}
\begin{abstract}
This note introduces the \textit{Macaulay2} package SchurVeronese, which gathers together data about Veronese syzygies and makes it readily accessible in \textit{Macaulay2}.   In addition to standard Betti tables, the package includes information about the Schur decompositions of the various spaces of syzygies.  The package also includes a number of functions useful for manipulating and studying this data. 
\end{abstract}
\maketitle

In~\cite{bruceErmanGoldsteinYang} the authors used a combination of high-throughput and high-performance computation and numerical techniques to compute the Betti tables of $\PP^2$ under the $d$-fold Veronese embedding, as well as the Betti tables of the pushforwards of line bundles $\cO_{\PP^2}(b)$ under that embedding, for a number of values of $b$ and $d$.  These computations resulted in new data, such as Betti tables, multigraded Betti numbers, and Schur Betti numbers.    (When $b=0$, most the cases had been previously computed in~\cite{wcdl}.)
This note introduces the \textit{SchurVeronese} package for \textit{Macaulay2}, which makes this data readily accessible via \textit{Macaulay2} for further experimentation and study.

\section{Veronese Syzygies}
Throughout this section we fix $n\in \mathbb{N}$ and let $S=\CC[x_0,x_1,\dots,x_n]$ be the polynomial ring with the standard grading.  The \defi{$d$th Veronese module of $S$ twisted by $b$} is
\[
S(b;d)\coloneqq \bigoplus_{i\in \ZZ} S_{di+b}.
\]
If $b=0$, then $S(0;d)$ is the Veronese subring of $S$, and if $b\ne 0$ then $S(b;d)$ is an $S(0;d)$-module.  Moreoever, if we set $R=\Sym(S_d)$ to be the symmetric algebra on $S_d$, then we may consider $S(b;d)$ as a graded $R$-module. Geometrically, if $b=0$ this corresponds to the homogenous coordinate ring of $\PP^n$ under the $d$-fold embedding $\PP^n \to \PP^{\binom{n+d}{d}-1}$, and for other $b$ it corresponds to the pushforward of $\cO_{\PP^n}(b)$ under the $d$-fold embedding.

Our interest is in studying the syzygies of $S(b;d)$.  See the introduction of \cite{bruceErmanGoldsteinYang} for background on Veronese syzygies including a summary of known results.  Throughout this paper, we set  
$
K_{p,q}(\PP^n, b;d):=\Tor^{R}_{p}(S(b;d),\CC)_{p+q},
$
which is isomorphic to the vector space of degree $p+q$ syzygies of $S(b;d)$ of homological degree $p$. Using the standard conventions for graded Betti numbers, the rank of the vector space $K_{p,q}$ corresponds to the Betti number $\beta_{p,p+q}$, and we write
$
\beta_{p,p+q}(S(b;d))
\coloneqq \dim \Tor^{R}_{p}(S(b;d),\CC)_{p+q} 
=\dim K_{p,q}(\PP^n, b;d).
$
Following the usual \textit{Macaulay2} notation, the \defi{Betti table} of $S(b;d)$ will be the table where $\beta_{p,p+q}(S(b;d))$ is placed in the $(p,q)$-spot. 

Outside of the case $n=1$, the Betti tables of $S(b;d)$ are unknown even for modest values of $d$.  There is not even a conjecture about what the Betti table of $S(b;d)$ should be for $n=2$ and $d\geq 7$.  

This package provides an array of computed data about $S(b;d)$ in the case $n=2$ and for $0\leq b < d \leq 8$ (though the data are incomplete for some of the larger values of $d$).  
While computing this data, including the Schur functor decompositions, took substantial time, the resulting data are concise and easy to work with in \textit{Macaulay2}.  
The bulk of this package thus consists of these output data, which is included as auxiliary files.  The functions provided in this package make this data accessible in a user-friendly way. 
Our hope is that this will allow those interested in Veronese syzygies to make headway on formulating conjectures and proving results in this area.  Moreover, as new cases of Veronese syzygies are computed, these can easily be incorporated into future versions of the package.

\section*{Acknowledgments}  We thank Claudiu Raicu and Gregory G.~Smith for useful conversations.

\section{An overview of the data}
When computing data for $S(b;d)$ we always work under the hypothesis that $0\leq b<d,$ as the  Betti table of $S(b;d)$ and $S(b+d;d)$ differ only by a vertical shift.  We have included data for the cases $n=1$ and $d\leq 10$, although this can also easily be computed using the Eagon-Northcott complex.  The main data are for the cases $n=2$ and $0\leq b < d\leq 8$.  In~\cite{bruceErmanGoldsteinYang}, we obtained full computations for $d\leq 6$; moreover since those algorithms worked in parallel with respect to multidegrees, we obtained incomplete data for some cases where $d=7,8$, and we have included that partial data in this package as well.

The algorithms in ~\cite{bruceErmanGoldsteinYang} are a mix of symbolic and numeric algebra.  Thus some entries in the data are not provably correct, while others are.  One can determine precisely when $K_{p,q}\ne 0$ by combining~\cite[Remark~6.5]{ein-lazarsfeld}, ~\cite[Theorem~2.2]{green-II}, and \cite[Theorem~2.c.6]{green-I}.  Our computation of a $K_{p,q}$-group (and all related data such as the Schur functor decomposition) will be provably correct if and only if $K_{p+1,q-1}$ and $K_{p-1,q+1}$ both vanish; in cases where this does not occur, the data for $K_{p,q}$ may have been computed numerically, and thus may not be provably correct.  For a longer discussion of potential numerical error issues, see \cite[\S5.2]{bruceErmanGoldsteinYang}.

\section{Total Betti Tables}
The Betti table for $S(b;d)$ can be called up using the \verb+totalBettiTally+ command.  For example, the Betti table of $S(2;4)$ when $n=2$ is produced below. 
\begin{verbatim}
i6 : totalBettiTally(4,2,0)

            0  1   2    3    4    5    6    7    8    9  10 11 12
o6 = total: 1 75 536 1947 4488 7095 7920 6237 3344 1089 175 24  3
         0: 1  .   .    .    .    .    .    .    .    .   .  .  .
         1: . 75 536 1947 4488 7095 7920 6237 3344 1089 120  .  .
         2: .  .   .    .    .    .    .    .    .    .  55 24  3

o6 : BettiTally
\end{verbatim}
Note that this is purely numeric: the package does not produce a minimal free resolution; the function simply returns the Betti numbers obtained by a previous computation.  The command \verb+totalBetti+ is similar, but expresses the Betti numbers simply as a hash table.

There is also a distinction between the indexing conventions.  When working with hash tables, we follow the more concise $K_{p,q}$ indexing conventions, instead of the $\beta_{p,p+q}$ indexing conventions used for Betti tallies.  Thus for instance, in the above example, the Betti number $\beta_{2,3}$ would correspond to key $(2,\{3\},3)$ in the Betti tally, but in the hash table \verb+totalBetti+ it corresponds to key $(2,1)$:
\begin{verbatim}
i4 :E =  totalBetti(4,2,0);

i5 : E#(2,1)
o5 = 536
\end{verbatim}

If one tries to call a Betti table outside of the acceptable range of $n,b,d$, we return a message explaining the problem.
\begin{verbatim}
i10 : totalBettiTally(4,3,0)
o10 = Need n = 1 or 2
\end{verbatim}

As noted above,  there were instances where we were able to partially compute Betti tables, for instance in the case of the $7$-uple embedding of $\PP^2$.  In those cases, we have recorded the entries that we know, and we mark the unknown entries with ``infinity''.  For example:
\begin{verbatim}
i14 : B = totalBetti(7,2,0);

i15 : B#(4,1)
o15 = 1031184

i16 : B#(20,1)
o16 = infinity
o16 : InfiniteNumber
\end{verbatim}
Thus, in this case, we see that $\dim K_{4,1}(\PP^2,2;7)= 1031184$, but we were unable to compute $\dim K_{20,1}(\PP^2,2;7)$.

\section{Schur Decomposition}

When $n=2$ and $d\geq 5$ the Betti tables of $S(b;d)$ are often unwieldy to work with, as they and their entries tend to be quite large. For example, the Betti table of $S(0;6)$ has 26 columns and many of the entries are on the order of $10^7$. 

A more concise way of recording the syzygies would be to take into account the symmetries coming from representation theory.  The natural linear action of $\GL_{n+1}(\CC)$ on $S$ induces an action on each vector space $K_{p,q}(\PP^n, b;d)$.  We can thus decompose this as a direct sum of Schur functors of total weight $d(p+q)+b$ i.e.
\[
K_{p,q}(\PP^n, b;d) = \bigoplus_{\substack{|\lambda|=d(p+q)+b}} \bS_{\lambda}(\CC^{n+1})^{\oplus m_{p,\lambda}(\PP^n,b;d)},
\]
with $m_{p,\lambda}(\PP^n,b;d)$ being the \defi{Schur Betti numbers} and $\bS_{\lambda}$ being the Schur functor corresponding to the partition $\lambda$~\cite[p. 76]{fultonHarris91}. The Schur Betti numbers can be accessed via the \verb+schurBetti+ command, which returns a hash table whose keys correspond to pairs $(p,q)$ for which $K_{p,q}(\PP^n, b;d)\neq0$, and whose values are lists corresponding to the Schur decomposition of this syzygy module. 

For example, let us consider $K_{2,1}(\PP^2, 0;4)$, which is a vector space of dimension $536$.  As a representation of $\GL_{3}(\CC)$, it turns out to be the sum of $9$ distinct Schur functors, each appearing with multiplicity $1$:
\[
K_{2,1}(\PP^2,0;4) =   \bS_{(9, 2,1)}\oplus \bS_{(8, 4,0)}\oplus \bS_{(8, 3,1)}\oplus \bS_{(7,5,0)}\oplus \bS_{(7, 4,1)}\oplus \bS_{(7, 3,2)}\oplus \bS_{(6,5,1)}\oplus \bS_{(6,4,2)}\oplus \bS_{(5, 4,1)}.
\]
\begin{verbatim}
i26 : (schurBetti(4,2,0))#(2,1)

o26 = {({9, 2, 1}, 1), ({8, 4, 0}, 1), ({8, 3, 1}, 1), ({7, 5, 0}, 1),
      ---------------------------------------------------------------------
 ({7, 4, 1}, 1), ({7, 3, 2}, 1), ({6, 5, 1}, 1), ({6, 4, 2}, 1), ({5, 4, 3}, 1)}

o8 : List
\end{verbatim}

From this, it is easy to compute statistics such as the number of representations and the number of distinct representations appearing in the Schur decomposition of $K_{p,q}(n,b;d)$. 
The \textit{SchurVeronese} package provides commands for these. For instance, in our example above we see that:

 \begin{verbatim} 
i11 : (numDistinctRepsBetti(4,2,0))#(2,1)
o11 = 9
\end{verbatim}
We can also display the number of representations appearing in each entry of the Betti table.  In the following example, the first table counts distinct Schur functors and the second counts the number of Schur functors with multiplicity.
\begin{verbatim}
i29 : makeBettiTally numDistinctRepsBetti(4,2,0)

             0 1 2  3  4  5  6  7  8  9 10 11 12
o29 = total: 1 2 9 17 23 23 26 25 21 13  3  1  1
          0: 1 . .  .  .  .  .  .  .  .  .  .  .
          1: . 2 9 17 23 23 26 25 21 13  1  .  .
          2: . . .  .  .  .  .  .  .  .  2  1  1
          
i30 : makeBettiTally numRepsBetti(4,2,0)

             0 1 2  3  4  5  6  7  8  9 10 11 12
o30 = total: 1 2 9 28 55 79 86 69 38 14  3  1  1
          0: 1 . .  .  .  .  .  .  .  .  .  .  .
          1: . 2 9 28 55 79 86 69 38 14  1  .  .
          2: . . .  .  .  .  .  .  .  .  2  1  1          
\end{verbatim}
Thus, $K_{4,1}(\PP^2,0;4)$ is the sum of $55$ irreducible representations, $23$ of which are distinct.

\section{Multigraded Betti Numbers}
One can also specialize the action of $\GL_{n+1}(\CC)$ to the torus action via $(\CC^*)^{n+1}$. This gives a decomposition of $K_{p,q}(\PP^n, b;d)$ into a sum of $\ZZ^{n+1}$-graded vector spaces of total weight $d(p+q)+b$. Specifically, writing $\CC(-\mathbf{a})$ for the vector space $\CC$ together with the $(\CC^*)^{n+1}$-action given by $(\lambda_0,\lambda_1,\dots,\lambda_n)\cdot \mu = \lambda_0^{a_0}\lambda_1^{a_1}\cdots \lambda_n^{a_n}\mu$ we have
\[
K_{p,q}(\PP^n, b;d) = \bigoplus_{\substack{\mathbf{a}\in \ZZ^{n+1} \\|\mathbf{a}|=d(p+q)+b}} \CC(-\mathbf{a})^{\oplus \beta_{p,\mathbf{a}}(\PP^n,b;d)}
\]
as $\ZZ^{n+1}$-graded vector spaces, or equivalently as $(\CC^*)^{n+1}$ representations. 

The \textit{SchurVeronese} package produces these multigraded Betti numbers for a number of examples  via the \verb+multiBetti+ command. As \verb+schurBetti+ does, this command returns a hash table whose keys correspond to pairs $(p,q)$ for which $K_{p,q}(\PP^n, b;d)\neq0$, and whose values are multigraded Hilbert polynomials encoding the multigraded decomposition of $K_{p,q}(n,b;d)$. More specifically the value of \verb+(multiBetti(d,n,b))#(p,q)+ is the polynomial
$
\sum_{\substack{\mathbf{a}\in \ZZ^{n+1} \\|\mathbf{a}|=d(p+q)+b}} \beta_{p,\ba}(n,b;d)\mathbf{t}^{\ba}.
$
where $\mathbf{t}^{\ba}$ denotes $t_0^{a_0}t_1^{a_1}\cdots t_{n}^{a_n}$. 

For example, $K_{12,2}(2,0;4)$ is the following $3$-dimensional $\ZZ^3$-graded vector space:
\[
K_{12,2}(2,0;4)\cong\CC(-(19,19,18))\oplus\CC(-(19,18,19))\oplus\CC(-(18,19,19)).
\]
The following code computes this, illustrating that the multigraded Hilbert function for $K_{12,2}(2,0;4)$ is $t_0^{19}t_1^{19}t_2^{18}+t_0^{19}t_1^{18}t_2^{19}+t_0^{18}t_1^{19}t_2^{19}$. 

\begin{verbatim}
i4 : (multiBetti(4,2,0))#(12,2)

      19 19 18    19 18 19    18 19 19
o4 = t  t  t   + t  t  t   + t  t  t
      0  1  2     0  1  2     0  1  2

o4 : QQ[t , t , t ]
         0   1   2
\end{verbatim}

\end{document}